\documentclass[onecolumn,10pt]{article}
\usepackage{amsfonts}
\usepackage{amsmath,amssymb}
\usepackage[singlespacing]{setspace}

\input{tcilatex}
\begin{document}

\date{}
\title{{\Large \textbf{Inverse source problem for a time-fractional heat
equation with generalized impedance boundary condition}}}
\author{Muhammed \c{C}\.{I}\c{C}EK$^{1}$ and Mansur I. ISMAILOV$^{2}$ \\
$^{1}$Department of Mathematics, Bursa Technical University, Bursa, Turkey \\
$^{2}$Department of Mathematics,Gebze Technical University, \\
41400 Gebze, Kocaeli, TURKEY\\
E-mail: muhammed.cicek@btu.edu.tr and mismailov@gtu.edu.tr\\
\ }
\maketitle

\bigskip

\qquad {\footnotesize \textbf{Abstract. The paper considers an inverse
source problem for a one-dimensional time-fractional heat equation with the
generalized impedance boundary condition. The inverse problem is the time
dependent source parameter identification together with the temperature
distrubution from the energy measurement. The well-posedness of the inverse
problem is shown by applying the Fourier expansion in terms of
eigenfunctions of a spectral problem which has the spectral parameter also
in the boundary condition and by using the results on Volterra type integral
equation with the kernel may have a diagonal singularity. }}


\section{Introduction}

\qquad

\qquad The mathematical analysis of inverse problem for the fractional
diffusion equation is extensively investigated in the last decade. The first
theoretical results for the inverse problem of finding the coefficients in
fractional diffusion equation are obtained in \cite{yamamoto
inv,zhang,SaYa,rundel,rundel2,kirane amc}. The mathematical literature for
such kind of inverse problems are rapidly growing but without being
exhaustive we refer only to \cite%
{yamamoto2,yamamoto3,Wu,kirane2,LiYaYa,Ismailov2,fenger,FuOlKi}.

The inverse problems for time fractional diffusion equations can be obtained
from the classical inverse diffusion problems by replacing the time
derivative with a fractional derivative. For such a non-classical
derivative, some standard methods for treating the inverse problems cannot
be applied. Such a difference implies that the inverse problems for
fractional diffusion equations should be more difficult. The difficulty
comes from the definition of the fractional-order derivatives, which is
essentially an integral with the kernel of weak singularity.

Because there is a wide mathematical literature for inverse problems for
time fractional heat equation we intend to create comprehensive lists of the
references which study inverse problems of finding source term and the
references appling some methods of spectral analysis. In the paper \cite%
{yamamoto inv} the uniqueness theorem is proved by using expansion in terms
of eigenfunctions of suitable Sturm-Liouville problem along the
Gelfand-Levitan theory. Similar eigenfunction expansion result along the
analytic continuation and Laplace transform is used in determination of
space-dependent source term in a fractional diffusion equation in \cite%
{zhang}. Spectral analysis of suitable Sturm- Liouville operator is actively
used in other coefficient identification problems for time-fractional
diffusion problems, \cite{rundel2,yamamoto3,Wu}. The papers \cite{kirane
amc, FuOlKi} and \cite{kirane2, Ismailov2} study inverse problems of finding
space dependent and time-dependent source terms, respectively, in
time-fractional diffusion equation by using eigenfunction expansion of the
non-self adjoint spectral problem along the generalized Fourier method.

We refer to \cite{zang,SaYa,rundel2,kirane
amc,Wu,kirane2,Ismailov2,fenger,FuOlKi} for the inverse problem of finding
time or space dependent source term associated with the time fractional heat
equation but this list is far from the complete. As for numerical methods
for such kind of inverse problems, see also \cite%
{murio1,murio2,twei,twei2,yamamoto4,BoIv}, and here we do not intend to
create any comprehensive lists of the references.

The rest of the paper is organized as follows: The mathematical formulation
of the problem is given in Section 2. In this section we also recall some
preliminary definitions and facts on fractional calculus, on linear
fractional differential equations and necessery properties of eigenvalues
and eigenfunctions of the auxiliary spectral problem. In Section 3, the
existence and uniqueness of the solution of the inverse time-dependent
source problem is proved. Finally, the continuous dependence upon the data
of the solution of the inverse problem is also shown in this section.

\section{Mathematical Preliminaries and Problem Formulation}

\subsection{Notes on Fractional Calculus\protect\bigskip}

\bigskip\ \textit{\qquad }In this part we recall some basic definitions and
facts on fractional calculus and present a necessary Lemma for further use.
For details see (\cite{ford, podlubny, samko}).

\textit{\qquad }Consider the following initial value problem , existence and
uniqueness result for such problem is given in \cite{ford} , for a linear
fractional differential equation with order $0<q<1,$ 
\begin{equation}
\left\{ 
\begin{array}{c}
D_{0+}^{q}[u(t)-u(0)]+\lambda u(t)=h(t),\text{ \ }t>0, \\ 
u(0)=u_{0},%
\end{array}%
\right.  \tag{2.1}
\end{equation}

\bigskip where $D_{0+}^{q}$ refers to the the Riemann-Liouville fractional
derivative of order $q$ ($0<q<1)$ in the time variable defined by%
\begin{equation*}
D_{0+}^{q}u(t)=\frac{1}{\Gamma (1-q)}\frac{d}{dt}\int_{0}^{t}\frac{u(\tau )}{%
(t-\tau )^{q}}d\tau \text{.}
\end{equation*}

The choice of the term $D_{0+}^{q}(u(x,.)-u(x,0))(t)$ instead of the usual
term $D_{0+}^{q}(u(x,.)(t)$ is not only to avoid the singularity at zero,
but also impose a meaningful initial condition (without fractional integral) 
\cite{ford}.

By using the Laplace transform, the solution of IVP (2.1) is given in \cite%
{Ismailov2} as

\begin{equation*}
u(t)=u_{0}E_{q,1}(-\lambda t^{q})+\int\limits_{0}^{t}(t-\tau )^{\alpha
-1}E_{q,q}(-\lambda (t-\tau )^{q})h(\tau )d\tau
\end{equation*}

\bigskip where $E_{q,\beta }(z)=\dsum\limits_{k=0}^{\infty }\frac{z^{k}}{%
\Gamma (qk+\beta )},\ \ q>0,\ \ \beta >0$ and $E_{q,1}(z)=E_{q}(z)=\dsum%
\limits_{k=0}^{\infty }\frac{z^{k}}{\Gamma (qk+1)},\ \ q>0$ are two
parameter and one parameter Mittag-Leffler function, respectively.

\bigskip \textit{\qquad }Let us introduce the functions $e_{q}(t,\lambda
):=E_{q}(-\lambda t^{q})$ and $e_{q,q}(t,\lambda ):=t^{q-1}E_{q,q}(-\lambda
t^{q})$ where $\lambda \in 
\mathbb{R}
_{+}$. Then following statements for the Mittag-Leffler type functions $%
e_{q}(t,\lambda )$ and $e_{q,q}(t,\lambda )$ holds.

\textbf{Proposition1 }(\cite{podlubny,mainardi})

$i)$ For $0<q<1$, $\lambda \in 
\mathbb{R}
_{+}$ the function $e_{q}(t,\lambda )$ is a monotonically decreasing
function.

$ii)$ The function $e_{q}(t,\lambda )$ has the estimates $e_{q}(t,\lambda
)\simeq e^{^{-\frac{\lambda }{\Gamma (q+1)}t^{q}}}$ $\ $for $t\ll 1$ and $%
e_{q}(t,\lambda )\simeq \frac{1}{\Gamma (1-q)\lambda t^{q}}$ for $t\gg 1.$

$iii)$ 
\begin{eqnarray*}
D_{0+}^{q}(e_{q,q}(t,\lambda )) &=&-\lambda e_{q,q}(t,\lambda ) \\
D_{0+}^{q}(e_{q}(t,\lambda )-e_{q}(0,\lambda )) &=&-\lambda e_{q}(t,\lambda )
\\
I_{0+}^{1-q}(e_{q,q}(t,\lambda )) &=&e_{q}(t,\lambda )
\end{eqnarray*}%
where $I_{0+}^{\gamma }\ $is the fractional integral of order $\gamma >0$
for an integrable function $f$ which is defined by $I_{0+}^{\gamma }f(t)=%
\frac{1}{\Gamma \left( \gamma \right) }\int\limits_{0}^{t}\left( t-s\right)
^{\gamma -1}f(s)ds.$

\bigskip Taking into also account the monotonically deceraesing character of 
$E_{q,q}(-\lambda t^{q})$ (\cite{s}), where $\lambda \in 
\mathbb{R}
_{+}$ the following statement holds true.

\textbf{Lemma 1}. (\cite{bzsy}) For $0<q<1$ Mittag-Leffler type function $%
E_{q,q}(-\lambda t^{q})$ satisfies 
\begin{equation*}
0\leq E_{q,q}(-\lambda t^{q})\leq \frac{1}{\Gamma (q)},\text{ \ \ \ \ \ }%
t\in \lbrack 0,\infty ),\text{ }\lambda \geq 0.
\end{equation*}

\bigskip

\bigskip \textit{\qquad }Now, we give the following lemma which is necessary
for further development.

\textbf{Lemma 2.} For $0<q<1$, $\lambda \in 
\mathbb{R}
_{+}$ we have 
\begin{equation*}
\int\limits_{t_{0}}^{t}(t-\tau )^{q-1}E_{q,q}(-\gamma (t-\tau )^{q})d\tau =%
\frac{1}{\gamma }(1-E_{q}(-\gamma (t-t_{0})^{q}).
\end{equation*}

\textbf{Proof. }By applying change of variable $z=t-\tau $ in the above
integral and using $\frac{d}{dt}E_{q}(-\gamma t^{q})=-\gamma
t^{q-1}E_{q,q}(-\gamma t^{q})$, \cite{samko} , we have 
\begin{eqnarray*}
\int\limits_{t_{0}}^{t}(t-\tau )^{q-1}E_{q,q}(-\gamma (t-\tau )^{q})d\tau
&=&\int\limits_{0}^{t-t_{0}}z^{q-1}E_{q,q}(-\gamma z^{q})dz \\
&=&-\frac{1}{\gamma }\int\limits_{0}^{t-t_{0}}\frac{d}{dz}E_{q}(-\gamma
z^{q})dz=\frac{1}{\gamma }(1-E_{q}(-\gamma (t-t_{0})^{q}).
\end{eqnarray*}

\bigskip

\bigskip \textit{\qquad }We give a rule for fractional differentiation with
order $0<q<1$ of an integral depending on a parameter, see \cite{podlubny},

\begin{equation*}
D_{0+}^{q}\int\limits_{0}^{t}K(t,\tau )d\tau =\int\limits_{0}^{t}D_{\tau
}^{q}K(t,\tau )d\tau +\underset{\tau \rightarrow t-0}{\lim }I_{\tau
}^{1-q}K(t,\tau ).\ 
\end{equation*}

We will also need to recall the following result.

\textbf{Lemma 3. }(\cite{samko})

\textit{\qquad }Let $f_{i}$ be a sequence of functions defined on the
interval $(a,b].$ Suppose the following conditions holds:

\bigskip $(i)$ $\ $the fractional derivative $D_{0+}^{q}$ $f_{i}(t)$ ,for a
given $q>0,$ exists for all $i\in 
\mathbb{N}
,t\in (a,b];$

$(ii)$ both series $\dsum \limits_{i=1}^{\infty }f_{i}(t)$ and $\dsum
\limits_{i=1}^{\infty }D_{0+}^{q}$ $f_{i}(t)$ are uniformly convergent on
the interval [$a+\epsilon ,b$] for any $\epsilon >0.$

\textit{\qquad }Then the function defined by the series $\dsum%
\limits_{i=1}^{\infty }f_{i}(t)$ is $q$ differentiable and satisfies $%
D_{0+}^{q}\dsum\limits_{i=1}^{\infty }f_{i}(t)=\dsum\limits_{i=1}^{\infty
}D_{0+}^{q}f_{i}(t).$

\subsection{Weak singular Volterra integral equations}

In this part we recall some basic results on Volterra type integral equation
with the kernel may have a diagonal singularity. For details see (\cite%
{vainikko,G}).

Consider the Volterra integral equation

\begin{equation}
u(t)=\int\limits_{0}^{t}Q(t,\tau )u(\tau )d\tau +f(t),\text{ }0\leq t\leq 1.
\tag{2.2}
\end{equation}%
Denote $\Delta =\{(t,\tau ):0\leq \tau <t\leq 1\}$ and introduce the class $%
S^{\nu }$of kernels $Q(t,\tau )$ that are defined and continuous on $\Delta $
and satisfy for $(t,\tau )\in \Delta $ the inequality%
\begin{equation*}
\left\vert Q(t,\tau )\right\vert \leq c\left( t-\tau \right) ^{-\nu },\text{ 
}\nu >0,\text{ }c=const>0.
\end{equation*}%
A kernel $Q(t,\tau )\in S^{\nu }$ is weakly singular if $\nu <1$. A weak
singularity of the kernel implies that the corresponding integral operator
is compact in the space $C[0,1]$. More precisely, the following statement
holds true.

\textbf{Lemma 4. }(\cite{vainikko}) Let $Q(t,\tau )\in S^{\nu }$ and $\nu <1$%
. Then the Volterra integral operator $T$ defined by $(Tr)(t)=\int%
\limits_{0}^{t}Q(t,\tau )r(\tau )d\tau $ maps $C[0,1]$ into itself and $%
T:C[0,1]\rightarrow C[0,1]$ is compact.

The proof of Lemma 2 is standard a detailed argument can be found in \cite%
{vainikko}. A consequence of Lemma 2 is the following result.

\bigskip

\textbf{Lemma 5. }(\cite{vainikko}) Let $f\in C[0,1]$ and $Q(t,\tau )\in
S^{\nu }$ with $\U{3bd} <1$. Then Eq. (2.2) has a unique solution $u\in
C[0,1]$.

We will also need to recall the results on weakly singular vesion of the
Gronwall's inequality.

\textbf{Lemma 6.} (\cite{G}) Let $T,\varepsilon ,M\in 
\mathbb{R}
_{+}$ and $0<q<1$. Moreover assume that $\delta :[0,T]\rightarrow 
\mathbb{R}
$ is a continuous function satisfying the inequality

\begin{equation*}
|\delta (t)|\leq \varepsilon +\frac{M}{\Gamma (q)}\int_{0}^{t}(t-\tau
)^{-\nu }|\delta (\tau )|d\tau ,\text{ with }\nu =1-q
\end{equation*}

for all $t\in \lbrack 0,T]$. Then 
\begin{equation*}
|\delta (t)|\leq \varepsilon E_{q}(Mt^{q})
\end{equation*}

for $t\in \lbrack 0,T]$.

\subsection{Problem Formulation and Associated Spectral Problem}

\qquad Let $T>0$ be a fixed number. In the rectangle $\Omega _{T}=\left\{
(x,t):0<x<1,0<t\leq T\right\} ,$ we will consider the following
time-dependent heat conduction equation 
\begin{equation}
D_{0+}^{q}(u(x,t)-u(x,0))=u_{xx}+r(t)f(x,t),\text{ \ \ \ \ }(x,t)\in \Omega
_{T},  \tag{2.3}
\end{equation}%
supplemented with the initial condition 
\begin{equation}
u(x,0)=\varphi (x),\text{ \ }0\leq x\leq 1,  \tag{2.4}
\end{equation}%
and the boundary condition 
\begin{equation}
u(0,t)=0,\text{ }0<t\leq T\text{\ ,}  \tag{2.5}
\end{equation}

\bigskip

\begin{equation}
au_{xx}(1,t)+du_{x}(1,t)+bu(1,t)=0\text{ , }0<t\leq T\text{\ ,}  \tag{2.6}
\end{equation}%
where $f(x,t),\varphi (x)$ are given functions and $a,b,d$ are given numbers
not simultaneously equal to zero.

The choice of the term $D_{0+}^{q}(u(x,.)-u(x,0))(t)$ instead of the usual
term $D_{0+}^{q}(u(x,.)(t)$ is not only to avoid the singularity at zero,
but also impose a meaningful initial condition (without fractional integral) 
\cite{ford}.

If the function $r(t)$ is known, the problem of finding $u(x,t)$ from
(2.3)-(2.6) is called the direct problem. However, the problem here is that
the source function $r(t)$ is unknown, which needs to be determined by
energy condition%
\begin{equation}
\int\limits_{0}^{1}u(x,t)dx=E(t),\text{ }0\leq t\leq T,  \tag{2.7}
\end{equation}%
where $E(t)$ are given functions. This problem is called the inverse
problem.\bigskip

On the other hand, use of integral condition (1.4) arises when the data on
the boundary cannot be measured directly, but only the average value of the
solution can be measured along the boundary. More precisely classical
boundary conditions (Neumann, Dirichlet and Robin type) are not always
adequate as it depends on the physical context which data can be measure at
the boundary of the physical domain. The classical boundary conditions
cases, one can have a selection of such large noise local space measurement,
but which on average produce a less noisy non-local measurement (2.7).

\qquad The inverse problem of finding $r(t)$ in classical heat conduction
equation $u_{t}=u_{xx}+r(t)f(x,t)$ with the conditions (2.3)-(2.7) was
already studied in \cite{Ismailov1} by using eigenfunction expansion in
terms of auxiliary spectral problem along the generalized Fourier method.
Our aim in present paper is to transfer this method to fractional cases.

Because the function $r$ is space independent, $a$, $b$ and $d$ are
constants and the boundary conditions are linear and homogeneous, the method
of separation of variables is suitable for studying the inverse problem
(2.3)-(2.7). It is well-known that the main difficulty in applying the
Fourier method is its basisness, i.e. the expansion in terms of
eigenfunctions of the auxiliary spectral problem \cite{kerimov}

\begin{equation}
\left\{ 
\begin{array}{l}
X^{\prime \prime }(x)+\mu X(x)=0,\text{ }0\leq x\leq 1, \\ 
\text{ }X(0)=0 \\ 
(a\mu -b)X(1)=dX^{^{\prime }}(1).%
\end{array}%
\right.  \tag{2.8}
\end{equation}%
In contrast to the classical Sturm--Liouville problem, this problem has the
spectral parameter also in the boundary condition. It makes it impossible to
apply the classical results on expansion in terms of eigenfunctions.

\qquad Consider the spectral problem (2.8) with $ad>0.$ It is known that its
eigenvalues $\mu _{n},$ $n=0,1,2...$ are real and simple. They form an
unbounded increasing sequence and the eigenfunction $X_{n}(x)$ corresponding
to $\mu _{n}$ has exactly $n$ simple zeros in the interval $(0,1)$. The sign
of the first eigenvalue $\mu _{0}$ is given as

\begin{eqnarray*}
\mu _{0} &<&0<\mu _{1}<\mu _{2}<...,\text{ if }-\frac{b}{d}>1, \\
\mu _{0} &=&0<\mu _{1}<\mu _{2}<...,\text{ if }-\frac{b}{d}=1, \\
0 &<&\mu _{0}<\mu _{1}<\mu _{2}<...,\text{ if }-\frac{b}{d}<1,
\end{eqnarray*}%
The eigenvalues and eigenfunctions have the following asymptotic behaviour: $%
\sqrt{\mu _{n}}=\pi n+O(\frac{1}{n}),$ $X_{n}(x)=\sin (\pi nx)+O(\frac{1}{n}%
),$ $\ $for sufficiently large $n$. Let $n_{0}$ be arbitrary fixed
non-negative integer. The system of eigenfunctions $\{X_{n}(x)\}$ $%
(n=0,1,2,...;n\neq n_{0})$ is a Riesz basis for $L_{2}[0,1].$ The system $%
\{u_{n}(x)\}$ $(n=0,1,2,...;n\neq n_{0})$ which is biorthogonal to the
system \ $\{X_{n}(x)\}$ $(n=0,1,2,...;n\neq n_{0})$\ \ has the form 
\begin{equation*}
u_{n}(x)=\frac{X_{n}(x)-\frac{X_{n}(1)}{X_{n_{0}}(1)}X_{n_{0}}(x)}{%
\left\Vert X_{n}\right\Vert _{L_{2}[0,1]}^{2}+\frac{a}{d}X_{n}^{2}(1)}.
\end{equation*}%
The following Bessel-type inequalities are valid for the systems $%
\{X_{n}(x)\}$ and $\{u_{n}(x)\}$ $(n=0,1,2,...;n\neq n_{0}).$

\bigskip

\bigskip \textbf{Lemma 7} \textit{(\cite{Ismailov1}) If }$\psi (x)\in
L_{2}[0,1]$\textit{, \ then the estimates }%
\begin{equation*}
\text{ }\dsum\limits_{\substack{ n=0  \\ n\neq n_{0}}}^{\infty }\left\vert
(\psi ,X_{n})\right\vert ^{2}\leq c_{1}\left\Vert \psi \right\Vert
_{L_{2}[0,1]}^{2},\text{ }\dsum\limits_{\substack{ n=0  \\ n\neq n_{0}}}%
^{\infty }\left\vert (\psi ,u_{n})\right\vert ^{2}\leq c_{2}\left\Vert \psi
\right\Vert _{L_{2}[0,1]}^{2}\text{ }
\end{equation*}%
\textit{hold for some positive constants }$c_{1}$ and $c_{2}$, \textit{where 
}$(\psi ,X_{n})=\int\limits_{0}^{1}\psi (x)X_{n}(x)dx$\textit{\ and }$(\psi
,u_{n})=\int\limits_{0}^{1}\psi (x)u_{n}(x)dx$ \textit{denote the usual
inner products in} $L_{2}[0,1].$

\bigskip Let us denote $\Phi _{n_{0}}^{4}[0,1]:=\{\psi (x)\in C^{4}[0,1];$ $%
\psi (0)=\psi ^{^{\prime \prime }}(0),$ $\psi (1)=\psi ^{^{\prime }}(1)=\psi
^{^{\prime \prime }}(1)=\psi ^{^{\prime \prime \prime }}(1)=0,$ $%
\int\limits_{0}^{1}\psi (x)X_{n_{0}}(x)dx=0\}.$

\textbf{Lemma 8} \textit{(\cite{Ismailov1})} \ \textit{If }$\psi (x)\in \Phi
_{n_{0}}^{4}[0,1]$\textit{, then we have:}%
\begin{eqnarray*}
\mu _{n}^{2}(\psi ,X_{n}) &=&(\psi ^{4},X_{n}),\text{ \ \ \ \ \ \ \ }\mu
_{n}^{2}(\psi ,u_{n})=(\psi ^{4},u_{n}),\text{ \ \ \ \ }n\geq 0,\text{ \ \ }
\\
\text{ }\dsum\limits_{\substack{ n=0  \\ n\neq n_{0}}}^{\infty }\left\vert
\mu _{n}(\psi ,X_{n})\right\vert &\leq &c_{3}\left\Vert \psi \right\Vert
_{C^{4}[0,1]},\text{ }\dsum\limits_{\substack{ n=0  \\ n\neq n_{0}}}^{\infty
}\left\vert \mu _{n}(\psi ,u_{n})\right\vert \leq c_{4}\left\Vert \psi
\right\Vert _{C^{4}[0,1]}, \\
\text{ }\dsum\limits_{\substack{ n=0  \\ n\neq n_{0}}}^{\infty }\left\vert
(\psi ,X_{n})\right\vert &\leq &c_{5}\left\Vert \psi \right\Vert
_{C^{4}[0,1]},\text{ }\dsum\limits_{\substack{ n=0  \\ n\neq n_{0}}}^{\infty
}\left\vert (\psi ,u_{n})\right\vert \leq c_{6}\left\Vert \psi \right\Vert
_{C^{4}[0,1]}\text{,}
\end{eqnarray*}%
where $c_{3},c_{4},c_{5}$ and $c_{6}$ \textit{are some positive constants}.

\bigskip

\section{Well-Posedness of the inverse problem}

\qquad The main result on existence and uniqueness of the solution of the
inverse problem (2.3)-(2.7) is presented. We prove the existence and
uniqueness of the solution of inverse problem (2.3)-(2.7) by means of
construction of a contraction mapping from energy condition (2.7). We call
classical solution as a pair of functions $\left\{ u(x,t),r(t)\right\} $
satisfying $u(x,t)\in C^{2}([0,1],%
\mathbb{R}
)$ , $D_{0+}^{q}(u(x,.)-u(x,0))\in C[(0,T],%
\mathbb{R}
]$ and $r(t)\in C[0,T]$.

\subsection{Formal construction of the solution}

\textbf{Theorem 1. }\textit{Suppose that the following conditions hold:}

$(A_{1})$\textit{\ }$\ \varphi (x)\in \Phi _{n_{0}}^{4}[0,1].$

\textit{\bigskip }$(A_{2})$\textit{\ }$\ E(t)\in C^{1}[0,T];$\textit{\ }$%
E(0)=\int\limits_{0}^{1}\varphi (x)dx.$

$(A_{3})$\textit{\ }$\ f(x,t)\in C(\overline{\Omega _{T}});$\textit{\ }$%
f(x,t)\in \Phi _{n_{0}}^{4}[0,1],\forall t\in $\textit{\ }$[0,T];$

\textit{\ \ \ \ \ \ }

\textit{\bigskip Then there exists a unique classical solution of the
inverse problem (2.3)-(2.7) in }$\Omega _{T}$\textit{.}

\textbf{Proof} \ For arbitrary $r(t)\in C[0,T]$, the solution of (2.3)-(2.6)
can be written in the form%
\begin{equation}
u(x,t)=\dsum\limits_{\substack{ n=0  \\ n\neq n_{0}}}^{\infty
}v_{n}(t)X_{n}(x),  \tag{3.1}
\end{equation}%
where the functions $v_{n}(t),$ $n=0,1,2,...;n\neq n_{0},$ are to be
determined. By using the Fourier's method, we can easily see that $v_{n}(t),$
$n=0,1,2,...;n\neq n_{0},$ satisfy the following system of countably many
linear fractional differential equations:%
\begin{equation}
\QATOPD\{ \} {D_{0+}^{q}(v_{n}(t)-v_{n}(0))+\mu
_{n}v_{n}(t)=r(t)f_{n}(t)}{v_{n}(0)=\varphi _{n}}  \tag{3.2}
\end{equation}%
where $\mu _{n}$ are eigenvalues of (2.8)

\qquad According to IVP (2.1), it can easily be seen that the solutions of
(3.2) is of the form 
\begin{equation}
v_{n}(t)=\varphi _{n}e_{q}(t,\mu _{n})+\int\limits_{0}^{t}e_{q,q}(t-\tau
,\mu _{n})r(\tau )f_{n}(\tau )d\tau  \notag
\end{equation}%
with $\varphi _{n}=\int\limits_{0}^{1}\varphi (x)u_{n}(x)dx,$ $%
f_{n}(t)=\int\limits_{0}^{1}f(x,t)u_{n}(x)dx$ and the solution is given
formally 
\begin{equation}
u(x,t)=\dsum\limits_{\substack{ n=0  \\ n\neq n_{0}}}^{\infty }[\varphi
_{n}e_{q}(t,\mu _{n})+\int\limits_{0}^{t}e_{q,q}(t-\tau ,\mu _{n})r(\tau
)f_{n}(\tau )d\tau ]X_{n}(x),  \tag{3.3}
\end{equation}

\bigskip \qquad The formulas (3.3) and (2.8) yield a following Volterra
integral equation of the first kind with respect to $r(t):$%
\begin{equation}
\int\limits_{0}^{t}K(t,\tau )r(\tau )d\tau +F(t)=E(t)  \tag{3.4}
\end{equation}%
where

\begin{equation*}
K(t,\tau )=\dsum\limits_{\substack{ n=0  \\ n\neq n_{0}}}^{\infty }\left[
f_{n}(\tau )e_{q,q}(t-\tau ,\mu _{n})\left(
\int\limits_{0}^{1}X_{n}(x)dx\right) \right]
\end{equation*}%
and 
\begin{equation*}
F(t)=\dsum\limits_{\substack{ n=0  \\ n\neq n_{0}}}^{\infty }\left[ \varphi
_{n}e_{q}(t,\mu _{n})\left( \int\limits_{0}^{1}X_{n}(x)dx\right) \right] .
\end{equation*}

$\bigskip $

Further, the equation (3.4) yields the following Volterra integral equation
of the second kind by taking fractional derivative $D_{0+}^{q}$:%
\begin{equation}
\int\limits_{0}^{t}D_{\tau }^{q}K(t,\tau )r(\tau )d\tau +r(t)\underset{\tau
\rightarrow t-0}{\lim }I_{\tau }^{q-1}K(t,\tau
)+D_{0+}^{q}(F(t)-F(0))=D_{0+}^{q}(E(t)-E(0)).  \tag{3.5}
\end{equation}%
\qquad \qquad

$\bigskip $By using the properties $(iii)$ in Proposition 1, it is easy to
show that 
\begin{eqnarray*}
\underset{\tau \rightarrow t-0}{\lim }I_{\tau }^{q-1}K(t,\tau ) &=&\underset{%
\tau \rightarrow t-0}{\lim }\dsum\limits_{\substack{ n=0  \\ n\neq n_{0}}}%
^{\infty }\left[ f_{n}(\tau )e_{q}(t-\tau ,\mu _{n})\left(
\int\limits_{0}^{1}X_{n}(x)dx\right) \right] \\
&=&\int\limits_{0}^{1}f(x,t)dx\text{,}
\end{eqnarray*}

\bigskip

\begin{equation}
D_{0+}^{q}(F(t)-F(0))=-\dsum\limits_{\substack{ n=0  \\ n\neq n_{0}}}%
^{\infty }\left[ \varphi _{n}\mu _{n}e_{q}(t,\mu _{n})\left(
\int\limits_{0}^{1}X_{n}(x)dx\right) \right]  \tag{3.6}
\end{equation}%
and 
\begin{equation}
D_{\tau }^{q}K(t,\tau )=-\dsum\limits_{\substack{ n=0  \\ n\neq n_{0}}}%
^{\infty }\left[ \mu _{n}e_{q,q}(t-\tau ,\mu _{n})f_{n}(\tau )\left(
\int\limits_{0}^{1}X_{n}(x)dx\right) \right]  \tag{3.7}
\end{equation}%
\qquad \qquad \qquad \qquad \qquad \qquad\ \qquad

\qquad We obtain the Volterra integral equation of the second kind with
respect to $r(t)$ in the form 
\begin{equation}
r(t)=P(t)+\int\limits_{0}^{t}Q(t,\tau )r(\tau )d\tau  \tag{3.8}
\end{equation}%
with the free term $P(t)=\frac{D_{0+}^{q}(E(t)-E(0))-D_{0+}^{q}(F(t)-F(0))}{%
\int\limits_{0}^{1}f(x,t)dx}$ and kernel $Q(t,\tau )=-\frac{D_{\tau
}^{q}K(t,\tau )}{\int\limits_{0}^{1}f(x,t)dx}$.

\bigskip

\bigskip According to the Lemma 1, we estimate the kernel of (3.8) in the
following form: 
\begin{equation*}
\left\vert Q(t,\tau )\right\vert =\frac{\left\vert D_{\tau }^{q}K(t,\tau
)\right\vert }{\left\vert \int\limits_{0}^{1}f(x,t)dx\right\vert }\leq \frac{%
C}{(t-\tau )^{\mu }}.
\end{equation*}%
where $C=\frac{Mc_{4}\underset{t\in \lbrack 0,T]}{\max }\left\Vert
f(.,t)\right\Vert _{C^{4}[0,1]}}{\Gamma (q)\underset{t\in \lbrack 0,T]}{\min 
}\mathit{\ }\left\vert \int\limits_{0}^{1}f(x,t)dx\right\vert }$ with $M\geq
\left\vert X_{n}(x)\right\vert $, $\forall n\in 
\mathbb{N}
,$ $\forall x\in \lbrack 0,1]$ and $\mu =1-q$.

Because the kernel $Q(t,\tau )$ belongs to the class $S^{\mu }$ with $0<\mu
<1$ the Volterra integral equation (3.8) is weakly singular. Then it has a
unique solution $r\in C[0,1]$ according by Lemma 4 and 5.

\subsection{Existence and Uniqueness of the solution of the inverse problem}

\qquad First, let us show that the solution of inverse problem (2.3)-(2.7)
is unique. Suppose that there were two solutions pair $(r,u)$ and $(a,v)$ of
the inverse problem (2.3)-(2.7). Then from form of solution (3.3) and (3.8),
we have 
\begin{equation}
u(x,t)-v(x,t)=\dsum\limits_{\substack{ n=0  \\ n\neq n_{0}}}^{\infty }\left(
\int\limits_{0}^{t}e_{q,q}(t-\tau ,\mu _{n})f_{n}(\tau )[r(\tau )-a(\tau
)]d\tau \right) X_{n}(x)  \tag{3.10}
\end{equation}

\bigskip and%
\begin{equation}
r(t)-a(t)=\int\limits_{0}^{t}Q(t,\tau )\left[ r(\tau )-a(\tau )\right] d\tau
.  \tag{3.11}
\end{equation}

$\bigskip $

$\bigskip $Then (3.11) yields $r=a$ . After inserting $r=a$ in (3.10), we
have $u=v.$

$\bigskip $\qquad So far we have proved the uniqueness of the solution of
the inverse problem. Because the solution $u(x,t)$ is formally given by the
series form (3.3), we need to show that the series corresponding to $u(x,t),$
$u_{x}(x,t),$ $u_{xx}(x,t)$ and $D_{0+}^{q}(u(x,.)-u(x,0))$ represent
continues functions. Under the assumptions $(A_{1})$-$(A_{3})$ and Lemma 7$,$
for all $(x,t)\in \overline{\Omega _{T}},$ the series corresponding to $%
u(x,t)$ is bounded above by the series 
\begin{equation}
\dsum\limits_{\substack{ n=0  \\ n\neq n_{0}}}^{\infty }\left[ \text{ }%
\left\vert \varphi _{n}\right\vert +\frac{1}{\mu _{n}}\underset{t\in \lbrack
0,T]}{\max }\left\vert r(t)\right\vert \underset{t\in \lbrack 0,T]}{\max }%
\left\vert f_{n}(t)\right\vert \right]  \tag{3.12}
\end{equation}%
The majorizing series (3.12) is convergent by using Lemmas 7 and 8. This
implies that by the Weierstrass M-test, the series (3.3) is uniformly
convergent in the rectangle $\overline{\Omega _{T}}$ and therefore, the
solution $u(x,t)$ is continuous in the rectangle $\overline{\Omega _{T}}$.

The majorizing series for $x$-partial derivative is

\bigskip

\begin{equation}
\text{ }\dsum\limits_{\substack{ n=0  \\ n\neq n_{0}}}^{\infty }\left[
\left\vert \varphi _{n}\right\vert \sqrt{\mu _{n}}+\frac{1}{\sqrt{\mu _{n}}}%
\underset{t\in \lbrack 0,T]}{\max }\left\vert r(t)\right\vert \underset{t\in
\lbrack 0,T]}{\max }\left\vert f_{n}(t)\right\vert \right]  \tag{3.13}
\end{equation}%
for $xx$-partial derivative is \bigskip

\begin{equation}
\text{ }\dsum\limits_{\substack{ n=0  \\ n\neq n_{0}}}^{\infty }\left[
\left\vert \varphi _{n}\right\vert \mu _{n}+\underset{t\in \lbrack 0,T]}{%
\max }\left\vert r(t)\right\vert \underset{t\in \lbrack 0,T]}{\max }%
\left\vert f_{n}(t)\right\vert \right] .  \tag{3.14}
\end{equation}%
It can easily be seen that the series (3.13)--(3.14) are convergent by
employing Lemma 7 and 8, Schwarz inequality. Hence by theWeierstrass M-test,
the series obtained for $x$-partial and $xx$-partial derivatives of (3.3) is
uniformly convergent in the rectangle $\overline{\Omega _{T}}.$ Therefore,
their sums $u_{x}(x,t)$ and $u_{xx}(x,t)$ are continuous in $\overline{%
\Omega _{T}}$.

\qquad Now it remains to show that $q$-fractional derivative of the series $%
u(x,t)-u(x,0)$ represents continuous function on $\Omega _{T}.$ We will show
that for any $\epsilon >0$ and $t\in \lbrack \epsilon ,T]$, the following
series 
\begin{equation*}
\dsum\limits_{\substack{ n=0  \\ n\neq n_{0}}}^{\infty
}[D_{0+}^{q}(v_{n}(t)-v_{n}(0))]X_{n}(x)
\end{equation*}%
corresponding to $q$-fractional derivative of the function $u(x,t)-u(x,0)$
is uniformly convergent. Before it we need to recall the Lemma 3. Now we can
see that equation (3.2) yields 
\begin{equation*}
D_{0+}^{q}(v_{n}(t)-v_{n}(0))=-\mu _{n}\varphi _{n}e_{q}(t,\mu _{n})-\mu
_{n}\int\limits_{0}^{t}e_{q,q}(t-\tau ,\mu _{n})r(\tau )f_{n}(\tau )d\tau
+r(t)f_{n}(t).
\end{equation*}

\bigskip We have the following estimates 
\begin{equation*}
\left\vert D_{0+}^{q}(v_{n}(t)-v_{n}(0))\right\vert \leq \left\vert \varphi
_{n}\right\vert \mu _{n}e_{q}(\varepsilon ,\mu _{n})+2\underset{t\in \lbrack
0,T]}{\max }\left\vert r(t)\right\vert \underset{t\in \lbrack 0,T]}{\max }%
\left\vert f_{n}(t)\right\vert
\end{equation*}

\bigskip and we obtain a majorant series as follow 
\begin{equation}
\dsum\limits_{\substack{ n=0  \\ n\neq n_{0}}}^{\infty }\left[ \left\vert
\varphi _{n}\right\vert \mu _{n}e^{-\frac{\mu _{n}}{\Gamma (1+q)}\varepsilon
^{q}}+2\underset{t\in \lbrack 0,T]}{\max }\left\vert r(t)\right\vert 
\underset{t\in \lbrack 0,T]}{\max }\left\vert f_{n}(t)\right\vert \right] . 
\tag{3.15}
\end{equation}

\bigskip Consequently $D_{0+}^{q}(u(x,t)-u(x,0))$ is uniformly convergent in
the rectangle $\Omega _{T}.$

\section{$\protect\bigskip $Lipschitz stability of the solution of the
inverse problem}

\bigskip \qquad The following result on continuously dependence on the data
of solution of the inverse problem (1.1)-(1.4) holds.

\bigskip \textbf{Theorem 2 \ }\textit{Let }$\Im $\textit{\ be the class of
triples in the form of }$\Phi =\left \{ f,\varphi ,E\right \} $\textit{\
which satisfy the assumptions }$(A_{1})-(A_{3})$\textit{\ of Theorem 1 and}

\begin{equation*}
0<N_{0}\leq \underset{t\in \lbrack 0,T]}{\min }\mathit{\ }\left\vert
\int\limits_{0}^{1}f(x,t)dx\right\vert ,\text{ }\left\Vert f\right\Vert
_{C^{4,0}(\overline{\Omega _{T}})}\leq N_{1},\text{ \ \ \ \ }\left\Vert
\varphi \right\Vert _{C^{4}[0,1]}\leq N_{2},\text{ \ \ \ \ \ \ }\left\Vert
E\right\Vert _{C^{1}[0,T]}\leq N_{3},\text{\ \ \ \ \ \ }
\end{equation*}%
\textit{for some positive constants }$N_{i},i=0,1,2,3$.

\textit{\bigskip \qquad Then the solution pair }$(r,u)$\textit{\ of the
inverse problem (2.3)-(2.7) depends continuously upon the data in }$\Im $.

\bigskip \textbf{Proof.} Let $\Phi =\left\{ f,\varphi ,E\right\} $ and $%
\widetilde{\Phi }=\left\{ \widetilde{f},\widetilde{\varphi },\widetilde{E}%
\right\} \in \Im $ be two sets of data. Let us denote $\left\Vert \Phi
\right\Vert =\left\Vert f\right\Vert _{C^{4,0}(\overline{\Omega _{T}}%
)}+\left\Vert \varphi \right\Vert _{C^{4}[0,1]}$\ $+$\ $\left\Vert
E\right\Vert _{C^{1}[0,T]}.$

Let $(r,u)$ and \ $(\widetilde{r},\widetilde{u})$ be the solutions of the
inverse problems (2.3)-(2.7) corresponding to the data \ $\Phi $ and $%
\widetilde{\Phi }$, respectively.

\textit{\qquad }According to (3.8) we have

\bigskip 
\begin{equation}
r(t)=P(t)+\int\limits_{0}^{t}Q(t,\tau )r(\tau )d\tau \   \tag{3.16}
\end{equation}%
\ and \ \ \ \ 
\begin{equation}
\widetilde{r}(t)=\widetilde{P}(t)+\int\limits_{0}^{t}\widetilde{Q}(t,\tau )%
\widetilde{r}(\tau )d\tau .  \tag{3.17}
\end{equation}

Firstly, from equations (3.5)-(3.7) and using \ $E(t)\in C^{1}[0,T]$ , Lemma
2, it is easy to compute that

\begin{equation}
\left\vert D_{0+}^{q}(F(t)-F(0))\right\vert \leq N_{4}\text{ , \ \ \ }%
\left\vert D_{0+}^{q}(E(t)-E(0))\right\vert \leq N_{5}  \tag{3.18}
\end{equation}%
\ \ \ 

\bigskip where $N_{4}=c_{4}N_{2}M,$ $N_{5}=\frac{T^{q}}{q\Gamma (1-q)}N_{3}$.

Let us estimate the difference $r-\widetilde{r}.$

\textit{\qquad }From (3.16) and (3.17) we obtain%
\begin{equation}
r(t)-\widetilde{r}(t)=P(t)-\widetilde{P}(t)+\int\limits_{0}^{t}\left[
Q(t,\tau )-\widetilde{Q}(t,\tau )\right] r(\tau )d\tau +\int\limits_{0}^{t}%
\widetilde{Q}(t,\tau )\left[ r(\tau )-\widetilde{r}(\tau )\right] d\tau . 
\tag{3.19}
\end{equation}

\bigskip Let $\epsilon _{1}=:\left\Vert P-\widetilde{P}\right\Vert
_{C([0,T])}+\frac{T^{q}}{q\Gamma (q)}\frac{2N_{1}Mc_{4}}{N_{0}^{2}}%
\left\Vert f-\widetilde{f}\right\Vert _{C^{4,0}(\overline{\Omega _{T}}%
)}\left\Vert r\right\Vert _{C([0,T])}.$

Then denoting $R(t)=:\left\vert r(t)-\widetilde{r}(t)\right\vert $, equation
(3.19) implies the inequality 
\begin{equation*}
R(t)\leq \epsilon _{1}+\frac{\epsilon _{2}}{\Gamma (q)}\int%
\limits_{0}^{t}(t-\tau )^{q-1}R(\tau )d\tau
\end{equation*}%
where $\epsilon _{2}=\frac{N_{1}Mc_{4}}{N_{0}}$Then, a weakly singular
Gronwall's inequality, see Lemma 6, implies that 
\begin{equation*}
R(t)\leq \epsilon _{1}E_{q}(\epsilon _{2}t^{q}),\text{ \ \ \ \ \ \ }t\in
\lbrack 0,T].\text{ }
\end{equation*}

\bigskip

Finally using (1) and 3.19 we obtain It follows from (3.8) that%
\begin{equation}
\left\Vert r-\widetilde{r}\right\Vert _{C([0,T])}\leq \epsilon
_{3}(\left\Vert P-\widetilde{P}\right\Vert _{C([0,T])}+\left\Vert
r\right\Vert _{C([0,T])}\frac{T^{q}}{q\Gamma (q)}\frac{N_{1}Mc_{4}}{N_{0}^{2}%
}\left\Vert f-\widetilde{f}\right\Vert _{C^{4,0}(\overline{\Omega _{T}})}) 
\tag{3.20}
\end{equation}%
where $\epsilon _{3}=E_{q}(\epsilon _{2}T^{q}).$ Also one can estimate that

\bigskip

\begin{equation}
\left\Vert P-\widetilde{P}\right\Vert _{C([0,T])}\leq M_{1}\left\Vert f-%
\widetilde{f}\right\Vert _{C^{4,0}(\overline{\Omega _{T}})}+M_{2}\left\Vert
\varphi -\widetilde{\varphi }\right\Vert _{C^{4}[0,1]}+M_{3}\left\Vert E-%
\widetilde{E}\right\Vert _{C^{1}([0,T])}  \tag{3.21}
\end{equation}%
\bigskip

where $M_{1}=\frac{N_{4}+N_{5}}{N_{0}^{2}},$ $M_{2}=\frac{c_{4}MN_{1}}{%
N_{0}^{2}},$ $M_{3}=\frac{T^{q}}{q\Gamma (1-q)}\frac{N_{1}}{N_{0}^{2}}.$

$\bigskip $By using the inequality (3.21), from (3.20) we obtain

$\bigskip \left\Vert r-\widetilde{r}\right\Vert _{C([0,T])}\leq M_{4}\left(
\left\Vert f-\widetilde{f}\right\Vert _{C^{4,0}(\overline{\Omega _{T}}%
)}+\left\Vert \varphi -\widetilde{\varphi }\right\Vert
_{C^{4}[0,1]}+\left\Vert E-\widetilde{E}\right\Vert _{C^{1}([0,T])}\right) $

\ \ \ \ \ \ \ \ \ \ \ \ \ \ \ \ \ \ \ \ \ \ \ \ \ \ $=M_{4}\left\Vert \Phi -%
\widetilde{\Phi }\right\Vert $

where $M_{4}=\max (\epsilon _{3}M_{2},\epsilon _{3}M_{3},\epsilon
_{3}M_{1}+\epsilon _{3}\left\Vert r\right\Vert _{C([0,T])}\frac{T^{q}}{%
q\Gamma (q)}\frac{N_{1}Mc_{4}}{N_{0}^{2}}).$

$\bigskip $This shows that $r$ depends continuously upon the input data.
From (3.3), a similar estimate is also obtained for the difference $u-%
\widetilde{u}$ as

\bigskip $\left\Vert u-\widetilde{u}\right\Vert _{C(\overline{\Omega _{T}}%
)}\leq M_{5}\left\Vert \Phi -\widetilde{\Phi }\right\Vert .$

This completes the proof of Theorem 2.

\bigskip

\section{Conclusion}

\bigskip The paper considers an inverse source problem of identification the
time dependent source parameter together with the temperature distrubution
from the energy measurement for a one-dimensional time-fractional heat
equation with the generalized impedance boundary condition. The
well-posedness of the inverse problem is shown by using the Fourier
expansion in terms of eigenfunctions of a spectral problem which has the
spectral parameter also in the boundary condition. To contrast to the fact
that the problem under consideration in present paper is obtained from the
classical diffusion problem \cite{Ismailov1} by replacing the time
derivative with a fractional derivative, the fractional inverse problem
should be more difficult. The difficulty comes from the definition of the
fractional-order derivatives, which is essentially an integral with the
kernel of weak singularity. For such a non-classical derivative, some
standard methods on Volterra type integral equations for treating the
inverse problems cannot be applied. It needs more difficult informations on
Volterra operators with the kernel of weak singularity at diagonal. This
approach can be extended to the time-fractional analogues of the classical
inverse initial boundary value problems for the heat equation with different
boundary and overdetermination conditions, which are the line of future
investigations.


\begin{thebibliography}{99}
\bibitem{yamamoto inv} Cheng J, Nakagawa J, Yamamoto M, Yamazaki T.
Uniqueness in an inverse problem for a one-dimensional fractional diffusion
equation. Inverse Problems 2009;25:115002.

\bibitem{yamamoto2} Li G, Zhang D, Jia X, Yamamoto M. Simultaneous inversion
for the space-dependent diffusion coefficient and the fractional order in
the time-fractional diffusion equation. Inverse Problems 2013; 29:065014.

\bibitem{zhang} Zhang Y, Xu X. Inverse source problem for a fractional
diffusion equation. Inverse Problems 2011 ;27:035010.

\bibitem{rundel} Jin B, Rundell W. An inverse problem for a one-dimensional
time-fractional diffusion problem. Inverse Problems 2012;28 :075010.

\bibitem{rundel2} Jin B, Rundell W. An inverse source Sturm-Liouville
problem with a fractional derivative. Journal of Computational Physics
2012;231:4954-4966.

\bibitem{yamamoto3} Luchko Y., Rundell W., Yamamoto M., Zuo L. Uniqueness
and reconstruction of an unknown semilinear term in a time-fractional
reaction-diffusion equation. Inverse Problems 2013;29 (6).

\bibitem{Wu} Wu B., Wu S. Existence and uniqueness of an inverse source
problem for a fractional integro-differential equation. Math. Meth. Appl.
Sci. 2014;37:1147-1158.

\bibitem{kirane amc} Kirane M, Salman AM. Determination of an unknown source
term and the temperature distribution for the linear heat equation involving
fractional derivative in time. Applied Mathematics and Computation
2011;218:63-170.

\bibitem{kirane2} Aleroev TS, Kirane M, Salman AM. Determination of a source
term for a time fractional diffusion equation with an integral type-over
determining condition. Electronic Journal of Differential Equations 2013;
No.270:pp 1-16.

\bibitem{murio1} Murio DA. Stable numerical solution of a
fractional-diffusion inverse heat conduction problem. Computers and
Mathematics with Applications 2007;53:1492--501.

\bibitem{murio2} Murio DA. Time fractional IHCP with Caputo fractional
derivatives. Comput Math Appl 2008;56:2371--81.

\bibitem{twei} Wei T, Zhang ZQ. Reconstruction of a time-dependent source
term in a time-fractional diffusion equation. Engineering Analysis with
Boundary Elements 2013;37:23--31.

\bibitem{twei2} Wang JG, Zhou YB, Wei T. Two regularization methods to
identify a space-dependent source for the time-fractional diffusion
equation. Applied Numerical Mathematics 2013; 68:39-57.

\bibitem{yamamoto4} Wang W., Yamamoto M., Han B. Numerical method in
reproducing kernel space for an inverse source problem for the fractional
diffusion equation. Inverse Problems 2013;29.

\bibitem{ford} Deithelm K, Ford NJ. Analysis of fractional differential
equations. J. Math. Anal. Appl. 2002;265:229-248.

\bibitem{bennani} Roussy G, Bennani A, Thiebaut J. Temperature runaway of
microwave irridiated materials. J. Appl. Phys. 1987;62:1167-1170.

\bibitem{cannon} Cannon JR, Esteva SP, van der Hoek J. A Galerkin procedure
for the diffusion equation subject to the specification of mass. Siam J.
Numer. Anal.1987;24:499-514.

\bibitem{Ismailov1} Hazanee A., Lesnic D., Ismailov MI, Kerimov NB. An
inverse time dependent source term for the heat equation with a
non-classical boundary condition. Applied Mathematical Modelling
2015;39:6258-6272.

\bibitem{Ismailov2} Ismailov MI, Cicek M. Inverse source problem for a
time-farctional diffusion equation with nonlocal boundary conditions.
Applied Mathematical Modelling 2016;40:4891-4899.

\bibitem{podlubny} Podlubny I. Fractional Differential Equations. Academic
Press: San Diego, 1999.

\bibitem{mainardi} Mainardi F. On some properties of the Mittag-Leffler
function $E_{\alpha }(-t^{\alpha })$, completely monotone for $t>0$ with $%
0<\alpha <1$. Discrete Contin. Dyn. Syst. Ser. B 2014;19 (7): 2267-2278.

\bibitem{kerimov} Kerimov NB, Ismailov MI. Direct and inverse problems for
the heat equation with a dynamic-type boundary condition, The IMA Journal of
Applied Mathematics, 2015, v. 80 (5), pp. 1519-1533.

\bibitem{samko} Samko G, Kilbas AA, Marichev OI. Fractional Integrals and
Derivatives: Theory and Applications. Gordon and Breach: Yverdon, 1993.

\bibitem{fenger} Feng P., Karimov ET. Inverse source problems for
time-fractional mixed parabolic-hyperbolic-type equations, J. Inverse
Ill-Posed Probl. 23 (2015), no. 4, 339-353.

\bibitem{FuOlKi} Furati KM., Iyiola OS., Kirane M. An inverse problem for a
generalized fractional diffusion, Applied Mathematics and Computation, 249
(2014) 24--31.

\bibitem{BoIv} Bondarenko AN and Ivaschenko DS. Numerical methods for
solving inverse problems for time fractional diffusion equation with
variable coefficient, J. Inverse Ill-Posed Probl. 17 (2009 ), 419--440

\bibitem{LiYaYa} Liu JJ, Yamamoto M and Yan L. On the uniqueness and
reconstruction for an inverse problem of the fractional diffusion process,
Appl. Numer. Math. 87 (2015 ), 1--9.

\bibitem{SaYa} Sakamoto K and Yamamoto M. Inverse source problem with a
final overdetermination for a fractional diffusion equation, Math. Control
Relat. Fields 1 (2011), 509--518.

\bibitem{vainikko} Pedas A, Vainikko G. Integral equations with diagonal and
boundary singularities of the kernel. Z Anal Anwend, 2006, 25(4), 487--516.

\bibitem{bzsy} Bai Z, Zhang S, Sun S, Yin C, Monotone iterative method for
fractional differential equations, Vol. 2016, No.06, 1-8.

\bibitem{s} Schneider WR, Completely monotone generalized Mittag-Leffler
functions, Expo. Math.,14 (1996), 3-16.

\bibitem{G} Dixon, J., McKee, S., Weakly singular discrete Gronwall
inequalities. Z. Angew. Math.Mech. 1986, v. 66, 535--544.
\end{thebibliography}
\end{document}